\documentclass [12pt] {amsart}
\usepackage{amsthm,amsmath,amssymb,amscd,graphics,enumerate}
\usepackage[all]{xy}
\usepackage{latexsym}
\usepackage{verbatim}
\usepackage{enumerate}
\usepackage{amsthm}
\usepackage{amsmath}
\usepackage{amscd}

\newcounter{theorem}[section]
\numberwithin{equation}{section}

\newtheorem{Thm}[theorem]{Theorem}

\theoremstyle{definition}
\newtheorem{conjecture}[theorem]{Conjecture}

\theoremstyle{remark}

\theoremstyle{remark}

\newtheorem{remark}[subsection]{Remark}

\numberwithin{equation}{section}
{\theoremstyle{remark}
 
 }

{\theoremstyle{definition}
 }

\newcommand{\Y}{\mathcal{Y}}

\newcommand{\B}{\mathcal{B}}

\newcommand{\D}{\mathcal{D}}

\def\<{\left\langle}
\def\>{\right\rangle}

\begin{document}

\title[Gerbe duality and relative Gromov-Witten theory]{On gerbe duality\\ and relative Gromov-Witten theory}

\author{Xiang Tang}
\address{Department of Mathematics\\ Washington University\\ St. Louis\\ MO 63130\\ USA}
\email{xtang@math.wustl.edu}

\author{Hsian-Hua Tseng}
\address{Department of Mathematics\\ Ohio State University\\ 100 Math Tower, 231 West 18th Ave.\\Columbus\\ OH 43210\\ USA}
\email{hhtseng@math.ohio-state.edu}

\date{\today}

\begin{abstract}
We formulate and study an extension of gerbe duality to relative Gromov-Witten theory.
\end{abstract}
\maketitle
\section{Introduction}

In this short note, we propose a conjecture about relative Gromov-Witten theory on a general noneffective Deligne-Mumford stack. 

Throughout this note, we work over $\mathbb{C}$. Let $G$ be a finite group, $\B$ a smooth (proper) Deligne-Mumford stack, and $$\pi: \Y\to \B$$ a $G$-gerbe. In \cite{hhps}, the {\em dual} of $\pi:\Y\to \B$ is a pair $(\widehat{\Y}, c)$ where $\widehat{\Y}$ is a disconnected stack with an \'etale map $$\widehat{\pi}:\widehat{\Y}=\coprod_i \widehat{\Y}_i\to \B$$ and $c$ is a $\mathbb{C}^*$-valued $2$-cocycle on $\widehat{\Y}$. 

We briefly recall the construction of $\widehat{\Y}$ and $c$, and refer the readers to \cite{tt1} for the detail. We focus on describing $\Y$ locally. A chart of $\B$ looks like a quotient of $\mathbb{C}^n$ by a finite group $Q$ acting linearly on $\mathbb{C}^n$. And a $G$-gerbe over $[\mathbb{C}^n/Q]$ can be described as follows. There is a group extension 
\[
1\longrightarrow G\longrightarrow H\longrightarrow Q\longrightarrow 1.
\]
The group $H$ acts on $\mathbb{C}^n$ via its homomorphism to $Q$. Locally a chart of $\Y$ over $\B$ looks like
\[
[\mathbb{C}^n/H]\longrightarrow [\mathbb{C}^n/Q]. 
\]

To construct the dual $\widehat{\Y}$, we consider the space $\widehat{G}$, the (finite) set of isomorphism classes of irreducible $G$-representations. As $G$ is a normal subgroup of $H$, $H$ acts on $G$ by conjugation, which naturally gives an $H$ action on $\widehat{G}$. Furthermore, as $G$ acts its by conjugation, the $G$ action on $\widehat{G}$ is trivial. Therefore, the quotient group $Q$ acts on $\widehat{G}$. The dual space $\widehat{\Y}$ locally looks like  $[(\widehat{G}\times \mathbb{C}^n)/Q]$. We notice that the dual $\widehat{\Y}$ has a canonical map to the quotient $\widehat{G}/Q$, and therefore is a disjoint union of stacks over $\widehat{G}/Q$.

The construction of $c$ is from the Clifford theory of induced representations \cite{cl}. Given an irreducible $G$-representation $\rho$ on $V_\rho$, we want to introduce an $H$ representation.  Let $[\rho]$ be the corresponding point in $\widehat{G}$. Recall that the finite group $Q$ acts on $\widehat{G}$. We denote $
Q_{[\rho]}$ to be the stabilizer group of $Q$ action on $[\widehat{G}]$ at the point $[\rho]$. For any $q\in Q_{[\rho]}$, define $q(\rho)$, a representation of $G$ on $V_\rho$, by 
\[
q(\rho)(g)=\rho\big(q^{-1}(g)\big).
\]
As $q$ fixes $[\rho]$ in $\widehat{G}$, $q(\rho)$, as a $G$-representation, is equivalent to $\rho$. Therefore, there is an intertwining operator $T_q^{\rho}$ on $V_\rho$ such that 
\[
T_q^\rho \rho =q(\rho) T_q^\rho. 
\]
In general, the operators $\{T_q^\rho\}_{q\in Q_{[\rho]}}$ fails to satisfies
\[
T_{q}^{\rho}\circ T^\rho_{\rho'}=T^{\rho}_{qq'}.
\]
By Schur's lemma, we can check that there is a number $c^{[\rho]}(q,q')\in \mathbb{C}^*$, such that 
\[
T_{q}^{\rho}\circ T^\rho_{\rho'}=c^{[\rho]}(q,q')T^{\rho}_{qq'}.
\]
In \cite{tt1}, we explained that this function $c^{[\rho]}(q,q')$ glues to a globally defined $\mathbb{C}^*$-gerbe over the dual $\widehat{\Y}$.

The authors of \cite{hhps} proposes a gerbe duality principle which suggests that there is an equivalence between the geometry of $\Y$ and that of the pair $(\widehat{\Y}, c)$. Several aspects of such an equivalence have been proven in \cite{tt1}. In \cite[Conjecture 1.8]{tt1}, the following conjecture is also explicitly formulated:

\begin{conjecture}\label{conj:gerbeGW}
As generating functions, the genus $g$ Gromov-Witten theory of $\Y$ is equal to the genus $g$ Gromov-Witten theory\footnote{It is also called $c$-twisted Gromov-Witten theory of $\widehat{\Y}$.} of $(\widehat{\Y}, c)$, $$GW_g(\Y)=GW_g(\widehat{\Y}, c).$$
\end{conjecture}
Conjecture \ref{conj:gerbeGW} has been proven in increasing generalities, see \cite{ajt1}, \cite{ajt2}, \cite{ajt3}, and \cite{tt2}. In particular, we have obtained the following theorem. 

\begin{Thm}(\cite[Theorem 1.1]{tt2})\label{thm:trivial-band}
When $\Y$ is 	a banded $G$-gerbe over $\B$, Conjecture \ref{conj:gerbeGW} holds true.
\end{Thm}

A toric Deligne-Mumford stack $\Y$ is a banded gerbe over an effective DM stack $\B$, see e.g. \cite{fmn}. As a corollary to Theorem \ref{thm:trivial-band}, we can compute the Gromov-Witten theory of $\Y$ in term of the (twisted) Gromov-Witten theory of the dual toric DM stack $\widehat{\Y}$.

From the perspective of Gromov-Witten theory, Gromov-Witten theory {\em relative} to a divisor is important. Therefore it is natural to consider an extension of Conjecture \ref{conj:gerbeGW} to the relative setting. Let $$D\subset \B$$ be a smooth (irreducible) divisor. The inverse images $$\D:=\pi^{-1}(D)\subset \Y, \quad \widehat{\D}:=\widehat{\pi}^{-1}(D)\subset \widehat{\Y}$$ are smooth divisors. The $c$-twisted Gromov-Witten theory of the relative pair $(\widehat{\Y}, \widehat{\D})$ can be defined using the construction of  \cite{af}, \cite{li1}, \cite{li2}, and \cite{pry}. A natural extension of Conjecture \ref{conj:gerbeGW} is the following  

\begin{conjecture}\label{conj:gerbeGWrel}
As generating functions, the genus $g$ Gromov-Witten theory of $(\Y, \D)$ is equal to the genus $g$ Gromov-Witten theory of $((\widehat{\Y},\widehat{\D}), c)$. Symbolically, $$GW_g(\Y,\D)=GW_g((\widehat{\Y}, \widehat{\D}), c).$$
\end{conjecture}

The purpose of this note is to present some evidence to Conjecture \ref{conj:gerbeGWrel}.

\section{Evidence to Conjecture \ref{conj:gerbeGWrel}}
By taking the divisor $D$ to be empty, we see that Conjecture \ref{conj:gerbeGWrel} implies Conjecture \ref{conj:gerbeGW}.
Next, we explain how to derive Conjecture \ref{conj:gerbeGWrel} from the full strength of Conjecture \ref{conj:gerbeGW}.

Let $r\geq 1$. We consider the stack of $r$-th roots of $\Y$ along $\D$, denoted by $$\Y_{\D, r}.$$ Consider also the stack of $r$-th roots of $\B$ along $D$, denoted by $$\B_{D,r}.$$ Let $\rho: \B_{D,r}\to B$ be the natural map. Consider the Cartesian diagram:
\begin{equation*}
\xymatrix{
 \rho^*\Y\ar[r]\ar[d]& \Y\ar[d]^\pi\\
\B_{D,r}\ar[r]^{\rho} & B.
}
\end{equation*}
The pull-back $\rho^*\Y\to \B_{D,r}$ is a $G$-gerbe. By functoriality property of root constructions, there is a natural map $$\Y_{\D,r}\to \rho^*\Y,$$ which is an isomorphism. Therefore we have 
\begin{equation}\label{eqn1}
GW_g(\Y_{D,r})=GW_g(\rho^*\Y).
\end{equation}

Since $\rho^*\Y\to \B_{D,r}$ is a $G$-gerbe, applying Conjecture \ref{conj:gerbeGW}, we have 
\begin{equation}\label{eqn2}
GW_g(\rho^*\Y)=GW_g(\widehat{\rho^*\Y}, c'),
\end{equation}
where $(\widehat{\rho^*\Y}, c')$ is the dual pair of the gerbe $\rho^*\Y$. By construction of the dual pair, we have $$\widehat{\rho^*\Y}=\rho^*\widehat{\Y}$$ and $c'=\rho^*c$. Here $\rho^*\widehat{\Y}$ fits in the Cartesian diagram
$$
\xymatrix{
\rho^*\widehat{\Y}\ar[r]\ar[d] & \widehat{\Y}\ar[d]^{\widehat{\pi}}\\
\B_{D,r}\ar[r]^{\rho} & B.
}
$$

By functoriality property of root constructions, we have $\rho^*\widehat{\Y}\simeq \widehat{\Y}_{\widehat{\D}, r}$. Therefore 
\begin{equation}\label{eqn3}
GW_g(\widehat{\rho^*\Y}, c')=GW_g(\widehat{\Y}_{\widehat{\D}, r}, \rho^*c).
\end{equation}
Combining (\ref{eqn1})--(\ref{eqn3}), we have 
\begin{equation}\label{eqn4}
GW_g(\Y_{D,r})=GW_g(\widehat{\Y}_{\widehat{\D}, r}, \rho^*c).
\end{equation} 
The left hand side of (\ref{eqn4}), $GW_g(\Y_{D,r})$, is a polynomial in $r$ for $r$ large. Furthermore, taking the $r^0$-coefficient, we obtain the relative Gromov-Witten invariants:
\begin{equation}\label{eqn5}
\text{Coeff}_{r^0}GW_g(\Y_{D,r})=GW_g(\Y, \D).
\end{equation}
This follows from the arguments of \cite{ty}, suitably extended to the setting of Deligne-Mumford stacks, see \cite{ty2}, \cite{cdw}.

The right hand side of (\ref{eqn4}), $GW_g(\widehat{\Y}_{\widehat{\D}, r}, \rho^*c)$, is also a polynomial in $r$ for $r$ large. Taking the $r^0$-coefficient, we obtain the relative Gromov-Witten invariants:
\begin{equation}\label{eqn6}
\text{Coeff}_{r^0}GW_g(\widehat{\Y}_{\widehat{\D}, r}, \rho^*c)=GW_g((\widehat{\Y}, \widehat{\D}), c).
\end{equation}
This again follows from a suitable extension of \cite{ty} as in \cite{ty2}, \cite{cdw}. Note that the $\rho^*c$-twist plays no role in applying the arguments of \cite{ty}, as they are essentially done at the level of virtual cycles (see \cite{fwy}), while $\rho^*c$-twist takes place at insertions.

Combining (\ref{eqn5}) and (\ref{eqn6}), we arrive at Conjecture \ref{conj:gerbeGWrel}.

\begin{remark}
In genus $0$, the argument about polynomiality in $r$ can be replaced by (a suitable extension of) the arguments of \cite{acw}.
\end{remark}

The following result follows directly from Theorem \ref{thm:trivial-band} and the above discussion. 
\begin{Thm}
When $\Y$ is a banded $G$-gerbe over $\B$, as generating functions, the genus $g$ Gromov-Witten theory of $(\Y, \D)$ is equal to the genus $g$ Gromov-Witten theory of $((\widehat{\Y},\widehat{\D}), c)$. Symbolically, $$GW_g(\Y,\D)=GW_g((\widehat{\Y}, \widehat{\D}), c).$$
\end{Thm}

\section{Acknowledgment}

X. Tang is supported in part by NSF Grant. H.-H. Tseng is supported in part by Simons Foundation Collaboration Grant.


\begin{thebibliography}{HLOY02}

\bibitem{acw} D. Abramovich, C. Cadman, J. Wise, {\em Relative and orbifold Gromov-Witten invariants}, Algebr. Geom. 4 (2017), no. 4, 472--500.

\bibitem{af} D. Abramovich, B. Fantechi, {\em Orbifold techniques in degeneration formulas}, Ann. Sc. Norm. Super. Pisa Cl. Sci. (5) 16 (2016), no. 2, 519--579.

\bibitem{ajt1} E. Andreini, Y. Jiang, H.-H. Tseng, {\em Gromov-Witten theory of root gerbes I: structure of genus 0 moduli spaces},  J. Differential Geom. Vol. 99, no. 1 (2015), 1--45.

\bibitem{ajt2} E. Andreini, Y. Jiang, H.-H. Tseng, {\em Gromov-Witten theory of product stacks}, Comm. Anal. Geom., Vol. 24 (2016), no. 2, 223--277.

\bibitem{ajt3} E. Andreini, Y. Jiang, H.-H. Tseng, {\em Gromov-Witten theory of banded gerbes over schemes}, arXiv:1101.5996.

\bibitem{cdw} B. Chen, C.-Y. Du, R. Wang, {\em Double ramification cycles with orbifold targets}, arXiv:2008.06484.

\bibitem{cl} A. Clifford, {\em Representations induced in an invariant subgroup}, Ann. of Math. (2) \textbf{38} (1937), no. 3, 533--550.

\bibitem{fwy} H. Fan, L. Wu, F. You, {\em Higher genus relative Gromov--Witten theory and DR-cycles}, J. London Math. Soc. 103 (2021), Issue 4, 1547--1576, arXiv:1907.07133.
  
\bibitem{fmn} B. Fantechi, E. Mann, F. Nironi, {\em Smooth toric Deligne-Mumford stacks}, J. Reine Angew. Math. 648 (2010), 201--244.

\bibitem{hhps} S. Hellerman, A. Henriques, T. Pantev, E. Sharpe, {\em Cluster decomposition, T-duality, and gerby CFTs}, Adv. Theor. Math. Phys., 11 (5) (2007), 751--818.

\bibitem{li1} J. Li, {\em Stable morphisms to singular schemes and relative stable morphisms}, J. Differential Geom. 57 (2001), 509--578.

\bibitem{li2} J. Li, {\em A degeneration formula of GW-invariants}, J. Differential Geom. 60 (2002), no. 2, 199--293.

\bibitem{pry} J. Pan, Y. Ruan, X. Yin, {\em Gerbes and twisted orbifold quantum cohomology}, Sci. China Ser. A, 51 (6) (2008), 995--1016.

\bibitem{tt1} X. Tang, H. -H. Tseng, {\em Duality theorems for \'etale gerbes on orbifolds}, Adv. Math. 250 (2014), 496--569.

\bibitem{tt2} X. Tang, H. -H. Tseng, {\em A quantum Leray-Hirsch theorem for banded gerbes}, J. Differential Geom. 119 (3), 459--511, 2021, arXiv:1602.03564. 

\bibitem{ty} H.-H. Tseng, F. You, {\em Higher genus relative and orbifold Gromov-Witten invariants}, Geom. Topol. 24 (2020) 2749--2779, arXiv:1806.11082.

\bibitem{ty2} H.-H. Tseng, F. You, {\em A Gromov-Witten theory for simple normal-crossing pairs without log geometry}, arXiv:2008.04844.

\end{thebibliography}
\end{document}